\documentclass[leqno]{amsart}
       
\usepackage{amsmath,amsfonts,amssymb}
\usepackage{amsthm}
\usepackage{hyperref}

\def\R{{\mathbb R}}
\def\N{{\mathbb N}}

\def\Sch{{\mathcal S}} 
\def\virgp{\raise 2pt\hbox{,}}

\def\bu{{\bf u}}

\def\({\left(}
\def\){\right)}
\def\<{\left\langle}
\def\>{\right\rangle}
\def\le{\leqslant}
\def\ge{\geqslant}

\def\Tend#1#2{\mathop{\longrightarrow}\limits_{#1\rightarrow#2}}

\def\d{{\partial}}

\def\eps{\varepsilon}

\def\si{{\sigma}}

\def\O{\mathcal O}

\DeclareMathOperator{\RE}{Re}
\DeclareMathOperator{\IM}{Im}

\theoremstyle{plain}
\newtheorem{theorem}{Theorem}[section]

\newtheorem{corollary}[theorem]{Corollary}
\newtheorem{proposition}[theorem]{Proposition}

\theoremstyle{remark}
\newtheorem{remark}[theorem]{Remark}

\numberwithin{equation}{section}

\begin{document}

\title[Instability for cubic NLS]{On the instability for the cubic
  nonlinear Schr\"odinger equation}
\author[R. Carles]{R{\'e}mi Carles}
\address{Universit\'e Montpellier~2\\ Math\'ematiques,
  UMR CNRS 5149\\ CC 051\\ 
   Place Eug\`ene Bataillon\\ 34095
  Montpellier cedex 5\\ France\footnote{Present address: Wolfgang
  Pauli Institute, Universit\"at Wien, 
        Nordbergstr.~15, A-1090 Wien}}
\email{Remi.Carles@math.cnrs.fr}
\thanks{Support by the ANR
  project SCASEN is acknowledged.} 
\begin{abstract}
We study the flow map associated to the cubic Schr\"odinger
equation in space dimension at least three. We consider initial data
of arbitrary size in $H^s$, where $0<s<s_c$, $s_c$ the critical index, and
perturbations in $H^\si$, where $\si<s_c$ is independent of $s$. We
show an instability mechanism in some Sobolev spaces of order smaller
than $s$. The analysis relies on two features of super-critical
geometric optics: creation of oscillation, and ghost effect. 
\end{abstract}
\subjclass[2000]{35B33; 35B65; 35Q55; 81Q05; 81Q20}
\maketitle

\section{Introduction}
\label{sec:intro}

We consider the Cauchy problem for the cubic,
defocusing Schr\"odinger equation: 
\begin{equation}
  \label{eq:nls}
  i\d_t \psi +\frac{1}{2}\Delta\psi = |\psi|^2\psi, \ x\in \R^n\quad
  ;\quad \psi_{\mid 
  t=0} =\varphi .
\end{equation}
Formally, the mass and energy associated to this equation are
independent of time:
\begin{align*}
&\text{Mass: }M[\psi](t) = \int_{\R^n}|\psi(t,x)|^2dx\equiv
M[\psi](0)= M[\varphi],\\
&\text{Energy: }E[\psi](t) = \int_{\R^n}|\nabla \psi(t,x)|^2dx
+\int_{\R^n}|\psi(t,x)|^{4}dx\equiv E[\psi](0)= E[\varphi].
\end{align*}
Scaling arguments yield the critical value for the Cauchy
problem in $H^s(\R^n)$:
\begin{equation*}
  s_c= \frac{n}{2}-1.
\end{equation*}
Assume $n\ge 3$, so that $s_c>0$. It was established in \cite{CW90}
that \eqref{eq:nls} is locally well-posed in $H^s(\R^n)$ if  $s\ge
s_c$. On the other hand, \eqref{eq:nls} is ill-posed in $H^s$ if
$s<s_c$ (\cite{CCT2}). Moreover, the following norm inflation
phenomenon was proved in \cite{CCT2} (see also \cite{BGTENS,CaARMA}):
if $0<s<s_c$, we can find $(\varphi_j)_{j 
  \in\N}$ in the 
Schwartz class $\Sch(\R^n)$ with 
\begin{equation}
  \label{eq:cipetite}
  \|\varphi_j\|_{H^s} \Tend j {+\infty} 0,
\end{equation}
 and a sequence of positive times $\tau_j\to 0$, such that the solution
$\psi_j$ to \eqref{eq:nls} with initial data
$\varphi_j$ satisfy:
\begin{equation*}
  \|  \psi_j(\tau_j) \|_{H^s} \Tend j {+\infty}
  +\infty.
\end{equation*}
In \cite{CaARMA}, this was improved to: we can find $t_j\to 0$ such that
\begin{equation*}
  \|  \psi_j(t_j) \|_{H^k} \Tend j {+\infty}
  +\infty, \quad \forall k \in \left] \frac{s}{1+s_c-s},s\right].
\end{equation*}
Note that \eqref{eq:cipetite} means that we consider the flow map near
the origin. We show that inside rings of $H^s$, the
situation is yet more involved: for data bounded in $H^s$, with
$0<s<s_c$, we consider 
perturbations which are small in $H^\si$ for \emph{any} $\si<s_c$, and
infer a similar conclusion. 
\begin{theorem}\label{theo:instab}
Let $n \ge 3$ and $0\le  s<s_c= \frac{n}{2}-1$. Fix $C_0,\delta>0$. We
can find two 
sequences of initial data $(\varphi_j)_{j \in\N}$ and $(\widetilde
\varphi_j)_{j \in\N}$ in the Schwartz class $\Sch(\R^n)$, with:
\begin{align*}
  &C_0-\delta \le \|\varphi_j\|_{H^s},\|\widetilde
  \varphi_j\|_{H^s}\le C_0+\delta \quad ;\quad
\|\varphi_j - \widetilde\varphi_j\|_{H^\si}\Tend j {+\infty} 0,\quad
  \forall \si <s_c, 
\end{align*}
and a sequence of positive times $t_j\to 0$, such that the solutions
$\psi_j$ and $\widetilde \psi_j$ to \eqref{eq:nls}, with initial data
$\varphi_j$ and $\widetilde \varphi_j$ respectively, satisfy:
\begin{equation*}
  \|  \psi_j(t_j)-\widetilde \psi_j(t_j) \|_{H^k} \Tend j {+\infty}
  +\infty,\quad \forall k \in \left] \frac{s}{1+s_c-s},s\right]\quad
  (\text{if }s>0),
\end{equation*}
\begin{equation*}
  \liminf_{j\to +\infty} \|  \psi_j(t_j)-\widetilde \psi_j(t_j)
  \|_{H^\frac{s}{1+s_c-s}} >0. 
\end{equation*}
\end{theorem}
The main novelty in this result is the fact that the initial data are
close to each other in $H^\si$, for any $\si<s_c$. In particular, this
range for $\si$ is independent of $s$. 
\begin{remark}
Like in \cite{BGTENS,CaARMA}, we consider initial data of the form 
  \begin{equation*}
    \varphi_j(x) = j^{\frac{n}{2}-s}a_0(jx),
  \end{equation*}
for some $a_0\in \Sch(\R^n)$ independent of $j$. The above
result holds for \emph{all} $a_0\in \Sch(\R^n)$ with,
say\footnote{Provided that we choose $j$  sufficiently large.},
$\|a_0\|_{H^s} =C_0$, and $\widetilde \varphi_j(x) =
(j^{\frac{n}{2}-s} + j)a_0(jx)$  (see Section~\ref{sec:reduc}).
\end{remark}
Considering the case $s=\frac{n}{4}$, we infer from the proof of
Theorem~\ref{theo:instab}: 
\begin{corollary}
 Let $n \ge 5$ and $C_0,\delta >0$. We can find two
sequences of initial data $(\varphi_j)_{j \in\N}$ and $(\widetilde
\varphi_j)_{j \in\N}$ in the Schwartz class $\Sch(\R^n)$, with:
\begin{align*}
C_0-\delta\le E[\varphi_j],E[\widetilde
 \varphi_j]\le C_0+\delta \quad ;\quad  M[\varphi_j]+M[\widetilde
 \varphi_j]+ E[\varphi_j - \widetilde \varphi_j] 
 \Tend j {+\infty} 0 ,  
\end{align*}
and a sequence of positive times $t_j\to 0$, such that the solutions
$\psi_j$ and $\widetilde \psi_j$ to \eqref{eq:nls} with initial data
$\varphi_j$ and $\widetilde \varphi_j$ respectively, satisfy:
\begin{equation*}
  \liminf_{j\to +\infty} E[\psi_j-\widetilde \psi_j](t_j)>0. 
\end{equation*} 
\end{corollary}

\section{Reduction of the problem: super-critical geometric optics}
\label{sec:reduc}
We now proceed as in \cite{CaARMA}. We set $\eps =
j^{s-s_c}$: $\eps\to 0$ as $j\to +\infty$. We change 
the unknown function as follows:
\begin{equation*}
  u^\eps (t,x) = j^{s-\frac{n}{2}} \psi_j \(
  \frac{t}{j^{s_c+2-s}}\virgp \frac{x}{j}\).  
\end{equation*}
Note that we have the relation:
\begin{equation*}
  \|\psi_j(t)\|_{\dot H^m} = j^{m-s}\left\|u^\eps\(
  j^{s_c+2-s}t \)\right\|_{\dot H^m}.
\end{equation*}
With initial data of the form
$\varphi_j(x)=j^{\frac{n}{2}-s}a_0(jx)+ j a_1(jx)$, \eqref{eq:nls} becomes:
\begin{equation}
  \label{eq:nlssemi}
  i\eps \d_t u^\eps +\frac{\eps^2}{2}\Delta u^\eps =
  |u^\eps|^2 u^\eps\quad ;\quad u^\eps(0,x) = a_0(x)+\eps a_1(x). 
\end{equation}
We emphasize two features for the WKB analysis associated to
\eqref{eq:nlssemi}. First, even if the initial datum is
independent of $\eps$, the solution instantly  becomes 
$\eps$-oscillatory. This is the argument of the proof of
\cite[Cor.~1.7]{CaARMA}. Second, the aspect which was not used in the
proof of \cite[Cor.~1.7]{CaARMA} is what was called \emph{ghost
  effect} in gas dynamics (\cite{Sone}): a perturbation of order
$\eps$ of the initial
datum may instantly become relevant at leading order. These two features
are direct 
consequences of the fact that \eqref{eq:nlssemi} is super-critical as
far as WKB analysis is concerned (see e.g. \cite{CaARMA}).
\smallbreak

Consider the two solutions $u^\eps$ and $\widetilde u^\eps$ of
\eqref{eq:nlssemi} with $a_1=0$ and $a_1=a_0$ respectively.
Then Theorem~\ref{theo:instab} stems from the following proposition,
which in turn is essentially a reformulation of \cite[Prop.~1.9
and 5.1]{CaARMA}. 
\begin{proposition}\label{prop:wnl}
  Let $n\ge 1$ and $a_0 \in \Sch (\R^n;\R)\setminus \{0\}$. There exist $T>0$
  independent of $\eps \in ]0,1]$, and 
  $a,\phi,\phi_1\in C([0,T];H^s)$ for all $s\ge 0$, such that: 
  \begin{equation*}
    \| u^\eps - a e^{i\phi/\eps}\|_{L^\infty([0,T];H^s_\eps)} + 
\| \widetilde u^\eps - ae^{i\phi_1}
e^{i\phi/\eps}\|_{L^\infty([0,T];H^s_\eps)} =\O(\eps),\quad \forall
s\ge 0,
  \end{equation*}
where 
\begin{equation*}
  \|f\|_{H^s_\eps}^2= \int_{\R^n} \( 1+ |\eps \xi|^2\)^{s}|\widehat
  f(\xi)|^2d\xi, 
\end{equation*}
and $\widehat f$ stands for the Fourier transform of $f$. In addition,
we have, in $H^s$:
\begin{equation*}
  \phi(t,x) = - t|a_0(x)|^2 + \O(t^3)\quad ;\quad \phi_1(t,x)
  =-2t|a_0(x)|^2 + \O(t^3) \quad \text{as }t\to 0.  
\end{equation*}
Therefore, there exists $\tau>0$ independent of $\eps$, such that:
\begin{equation*}
  \liminf_{\eps \to 0}\eps^{s}\| u^\eps(\tau) - \widetilde
  u^\eps(\tau)\|_{\dot H^s}>0, \quad \forall s\ge 0.
\end{equation*}
\end{proposition}

\section{Outline of the proof of Proposition~\ref{prop:wnl}}
\label{sec:outline}

The idea, due to E.~Grenier \cite{Grenier98}, consists in writing the
solution to 
\eqref{eq:nlssemi} 
as $u^\eps(t,x)= a^\eps(t,x)e^{i\phi^\eps(t,x)/\eps}$, where $a^\eps$
is complex-valued, and $\phi^\eps$ is real-valued. We assume that
$a_0, a_1\in \Sch(\R^n)$ are independent of $\eps$. For simplicity, we
also assume that they are real-valued. Impose:
\begin{equation}
  \label{eq:emmanuel}
  \left\{
    \begin{aligned}
      &\d_t \phi^\eps + \frac{1}{2}|\nabla \phi^\eps|^2
      +|a^\eps|^2=0\quad ;\quad \phi^\eps(0,x)=0.\\
& \d_t a^\eps +\nabla \phi^\eps\cdot \nabla a^\eps+\frac{1}{2}a^\eps
      \Delta \phi^\eps =i\frac{\eps}{2}\Delta a^\eps\quad ;\quad
      a^\eps(0,x)=a_0(x)+\eps a_1(x).  
    \end{aligned}
\right.
\end{equation}
Working with the unknown function $\bu^\eps = \,^t(\RE a^\eps, \IM
a^\eps, \d_1 \phi^\eps,\ldots,\d_n \phi^\eps)$, \eqref{eq:emmanuel}
yields a symmetric quasi-linear hyperbolic system: for $s>n/2+2$, there
exists $T>0$ independent of $\eps\in ]0,1]$ (and of $s$, from tame
estimates), such that \eqref{eq:emmanuel} has a unique solution
$(\phi^\eps,a^\eps)\in C([0,T];H^s)^2$. Moreover, the bounds in
$H^s(\R^n)$ are independent of $\eps$, and  we see that
$(\phi^\eps,a^\eps)$ converges to $(\phi,a)$, 
solution of:
\begin{equation}
  \label{eq:euler}
  \left\{
    \begin{aligned}
      &\d_t \phi + \frac{1}{2}|\nabla \phi|^2
      +|a|^2=0\quad ;\quad \phi(0,x)=0.\\
& \d_t a +\nabla \phi\cdot \nabla a+\frac{1}{2}a
      \Delta \phi =0\quad ;\quad
      a(0,x)=a_0(x).  
    \end{aligned}
\right.
\end{equation}
More precisely, energy estimates for symmetric systems yield: 
\begin{equation*}
  \|\phi^\eps -\phi\|_{L^\infty([0,T];H^s)}+ \|a^\eps
  -a\|_{L^\infty([0,T];H^s)} =\O(\eps),\quad \forall s\ge 0.  
\end{equation*}
One can prove that $\phi^\eps$ and $a^\eps$ have an asymptotic
expansion in powers of $\eps$. Consider the next term, given by:
\begin{equation*}
\left\{
  \begin{aligned}
    &\partial_t \phi^{(1)} +\nabla \phi \cdot \nabla \phi^{(1)} +
    2\operatorname{Re}\left(\overline a a^{(1)}\right)=0\  ;\ 
    \phi^{(1)}\big|_{t=0}=0.\\ 
   &\partial_t a^{(1)} +\nabla\phi\cdot \nabla a^{(1)} + \nabla
   \phi^{(1)}\cdot \nabla a + \frac{1}{2} a^{(1)}\Delta \phi
   +\frac{1}{2}a\Delta \phi^{(1)}      = \frac{i}{2}\Delta a\ 
    ;\  a^{(1)}\big|_{t=0}=a_1.
  \end{aligned}
\right.
\end{equation*}
Then $a^{(1)},\phi^{(1)}\in
L^\infty([0,T];H^s)$ for every $s\ge 0$, and
\begin{equation*}
  \|a^\eps - a - \eps a^{(1)}\|_{L^\infty([0,T_*];H^s)}+
  \|\Phi^\eps - \phi - \eps 
  \phi^{(1)}\|_{L^\infty([0,T_*];H^s)} \le C_s\eps^2,\quad
  \forall s\ge 0\, . 
\end{equation*}
Observe that since $a$ is real-valued, $(\phi^{(1)},\RE (\overline a
a^{(1)}))$ solves an 
homogeneous linear system. Therefore, if $\RE (\overline a
a^{(1)})=0$ at time $t=0$, then $\phi^{(1)}\equiv 0$.
\smallbreak

Considering the cases $a_1=0$ and $a_1=a_0$ for $u^\eps$ and
$\widetilde u^\eps$ respectively, we obtain the first assertion of
Prop.~\ref{prop:wnl}. Note that the above $\O(\eps^2)$ becomes
an $\O(\eps)$ only, since we divide $\phi^\eps$ and $\phi$ by $\eps$. 
This also explains why the first estimate of
Prop.~\ref{prop:wnl} is stated in $H^s_\eps$ and not in $H^s$. 
The rest of the proposition follows easily.  
\begin{remark}
  We could use the \emph{ghost effect} at higher order. For $N\in \N$,
  assume $\widetilde u^\eps_{\mid t=0}=(1+\eps^N)a_0$ for
  instance. Then for some $\tau>0$ independent of $\eps$, we have
\begin{equation*}
  \liminf_{\eps \to 0}\(\eps^{s}\| u^\eps(\tau) - \widetilde
  u^\eps(\tau)\|_{\dot H^s}\times \eps^{1-N}\)>0, \quad \forall s\ge 0.
\end{equation*}
Back to the functions $\psi$, the range for $k$ becomes:
\begin{equation*}
  k\ge \frac{s+(s_c-s)(N-1)}{1+s_c-s}\cdot
\end{equation*}
For this lower bound to be strictly smaller than $s$, we
have to assume $s> N-1$. 
\end{remark}


\begin{thebibliography}{1}

\bibitem{BGTENS}
N.~Burq, P.~G\'erard, and N.~Tzvetkov, \emph{Multilinear eigenfunction
  estimates and global existence for the three dimensional nonlinear
  {S}chr\"odinger equations}, Ann. Sci. \'Ecole Norm. Sup. (4) \textbf{38}
  (2005), no.~2, 255--301.

\bibitem{CaARMA}
R.~Carles, \emph{Geometric optics and instability for semi-classical
  {S}chr\"odinger equations}, Arch. Ration. Mech. Anal. (2007), to appear ({\tt
  doi:10.1007/s00205-006-0017-5}).

\bibitem{CW90}
T.~Cazenave and F.~Weissler, \emph{The {C}auchy problem for the critical
  nonlinear {S}chr\"odinger equation in ${H}^s$}, Nonlinear Anal. TMA
  \textbf{14} (1990), 807--836.

\bibitem{CCT2}
M.~Christ, J.~Colliander, and T.~Tao, \emph{Ill-posedness for nonlinear
  {S}chr\"odinger and wave equations}, Ann. Inst. H. Poincar\'e Anal. Non
  Lin\'eaire, to appear. See also {\tt arXiv:math.AP/0311048}.

\bibitem{Grenier98}
E.~Grenier, \emph{Semiclassical limit of the nonlinear {S}chr\"odinger equation
  in small time}, Proc. Amer. Math. Soc. \textbf{126} (1998), no.~2, 523--530.

\bibitem{Sone}
Y.~Sone, K.~Aoki, S.~Takata, H.~Sugimoto, and A.~V. Bobylev,
  \emph{Inappropriateness of the heat-conduction equation for description of a
  temperature field of a stationary gas in the continuum limit: examination by
  asymptotic analysis and numerical computation of the {B}oltzmann equation},
  Phys. Fluids \textbf{8} (1996), no.~2, 628--638.

\end{thebibliography}
\end{document}